\documentclass[12pt,a4paper,fleqn]{article}
\usepackage{amsfonts,amsmath,latexsym,amssymb,graphicx,mathrsfs}

\usepackage{graphicx}
\usepackage{color}
\usepackage{amssymb}
\usepackage{amssymb}
\usepackage[T1]{fontenc}
\usepackage{latexsym}
\usepackage{xypic}
\usepackage{eufrak}
\usepackage{euscript}
\usepackage{amsfonts,amsmath}
\usepackage{verbatim}
\usepackage[english]{babel}
\usepackage{mathrsfs}

\newtheorem{prop}{Proposition}[section]

\newtheorem{cor}[prop]{Corollary}
\newtheorem{lemma}[prop]{Lemma}
\newtheorem{rem}[prop]{Remark}

\newtheorem{thm}[prop]{Theorem}

\renewcommand{\geq}{\geqslant}
\def\leq{\leqslant}

\newcommand{\Z}{\mathbb{Z}}
\newcommand{\R}{\mathbb{R}}

\newcommand{\C}{\mathbb{C}}

\def\HH{\EuFrak H}

\def\e{\varepsilon}

\def\1{{\mathbf{1}}}

\def\1{{\mathbf{1}}}
\def\0.5{{\frac{1}{2}}}

\newcommand{\beas}{\begin{eqnarray*}}
\newcommand{\enas}{\end{eqnarray*}}

\newcommand{\NN}{\mathscr{N}}

\newcommand{\fin}
{ \vspace{-0.6cm}
\begin{flushright}
\mbox{$\Box$}
\end{flushright}
\noindent }

\newcommand{\qed}{\nopagebreak\hspace*{\fill}
{\vrule width6pt height6ptdepth0pt}\par}

\begin{document}

\begin{center}
{\Large{\bf Almost sure central limit theorems\\ on the Wiener space}}
\\~\\ by
Bernard Bercu\footnote{Institut de Math\'ematiques de Bordeaux, Universit\'e Bordeaux 1, 351 cours de
la lib\'eration, 33405 Talence cedex,
France. Email: \texttt{Bernard.Bercu@math.u-bordeaux1.fr}},
Ivan Nourdin\footnote{Laboratoire de Probabilit\'es et Mod\`eles Al\'eatoires, Universit\'e
 Pierre et Marie Curie (Paris VI), Bo\^ite courrier 188, 4 place Jussieu, 75252 Paris Cedex 05,
France. Email: \texttt{ivan.nourdin@upmc.fr}}
and Murad S. Taqqu\footnote{Boston University, Departement of Mathematics, 111 Cummington Road, Boston (MA), USA.
Email: \texttt{murad@math.bu.edu}}\footnote{Murad S. Taqqu was partially supported by the NSF Grant DMS-0706786 at
Boston University.}\\
{\small {\it Universit\'e Bordeaux 1, Universit\'e Paris 6 and Boston University}}\\~\\
\end{center}

{\small \noindent {\bf Abstract:}
In this paper, we study almost sure central limit theorems for sequences of functionals of general
Gaussian fields. We apply our result to non-linear functions of stationary Gaussian sequences.
We obtain almost sure central limit
theorems for these non-linear functions when they converge in law to a normal distribution.\\
}

{\small \noindent {\bf Key words}: Almost sure limit theorem; multiple stochastic integrals;
fractional Brownian motion; Hermite power variation.
\\


\noindent {\it This version}: December 11, 2009}

\section{Introduction}
Let $\{X_{n}\}_{n\geq 1}$ be a sequence of real-valued independent identically distributed random variables with
$E[X_{n}]=0$ and $E[X_{n}^{2}]=1$, and denote
\[
S_{n}= \frac1{\sqrt{n}}\sum_{k=1}^n  X_k.
\]
The celebrated almost sure central limit theorem (ASCLT)
states that the sequence of random empirical measures, given by
\[
\frac{1}{\log n}\sum_{k=1}^{n}\frac{1}{k}\delta_{S_k}
\]
converges almost surely to the $\NN(0,1)$ distribution as $n\to\infty$.
In other words, if $N$ is a $\NN(0,1)$ random variable, then, almost surely, for all $x\in\R$,
\[
\frac{1}{\log n}\sum_{k=1}^{n}\frac{1}{k}{\bf 1}_{\{S_k\leq x\}} \longrightarrow
P(N\leq x),\quad\mbox{as $n\to\infty$},
\]
or, equivalently, almost surely, for any bounded and continuous function $\varphi:\R\to\R$,
\begin{equation}
\label{ASCLTh}
\frac{1}{\log n}\sum_{k=1}^{n}\frac{1}{k}\varphi(S_{k}) \longrightarrow
E[\varphi(N)],\quad\mbox{as $n\to\infty$}.
\end{equation}

The ASCLT was stated first by L\'evy \cite{Levy} without proof.
It was then forgotten for half century.
It was rediscovered by Brosamler \cite{Brosamler} and Schatte \cite{Schatte}
and proven, in its present form, by Lacey and Philipp \cite{LaceyPhillip}.
We refer the reader to Berkes and Cs\'aki \cite{BerkesCsaki} for a universal
ASCLT covering a large class of limit theorems for partial sums, extremes, empirical
distribution functions and local times associated
with independent random variables $\{X_n\}$, as well as to the work of
Gonchigdanzan \cite{these}, where extensions of the ASCLT to weakly dependent
random variables are studied, for example in the context of strong mixing or $\rho$-mixing.
Ibragimov and Lifshits \cite{IL,IbragimovLifshits2} have provided a criterion
for (\ref{ASCLTh}) which does not require the sequence $\{X_n\}$ of random variables
to be necessarily independent nor
the sequence $\{S_n\}$ to take the specific form of partial sums.
This criterion is stated in Theorem \ref{thm-IL} below.

Our goal in the present paper is to investigate the ASCLT for a sequence of functionals of general
Gaussian fields. Conditions ensuring
the convergence in law of this sequence to the standard $\NN(0,1)$ distribution
may be found in \cite{NP07,NPR} by Nourdin, Peccati and Reinert.
Here, we shall propose a suitable criterion for
this sequence of functionals to satisfy also the ASCLT.
As an application, we shall consider some non-linear
functions of strongly dependent Gaussian random variables.

The paper is organized as follows. In Section 2,
we present the basic elements of Gaussian analysis and Malliavin calculus
used in this paper. An abstract version of our ASCLT is stated and proven in Section 3, as well as an application
to partial sums of non-linear functions of a strongly dependent Gaussian sequence.
In Section 4, we apply our ASCLT to discrete-time fractional Brownian motion.
In Section 5, we consider applications
to partial sums of Hermite polynomials of strongly dependent Gaussian sequences, when the limit
in distribution is Gaussian. Finally, in Section 6, we discuss the case where the
limit in distribution is non-Gaussian.

\section{Elements of Malliavin calculus}\label{sec22}
We shall now present the basic elements of Gaussian analysis and Malliavin calculus
that are used in this paper. The reader is referred to the monograph by  Nualart \cite{Nbook}
for any unexplained definition or result.

Let $\EuFrak H$ be a real separable Hilbert space. For any $q\geq 1$,
let $\EuFrak H^{\otimes q}$ be the $q$th tensor product of $\EuFrak H$ and denote
by $\EuFrak H^{\odot q}$ the associated $q$th symmetric tensor product.
We write $X=\{X(h),h\in \EuFrak H\}$ to indicate an isonormal Gaussian process over
$\EuFrak H$, defined on some probability space $(\Omega ,\mathcal{F},P)$.
This means that $X$ is a centered Gaussian family,
whose covariance is given in terms of the
scalar product of $\EuFrak H$ by
$E\left[ X(h)X(g)\right] =\langle h,g\rangle _{\EuFrak H}$.

For every $q\geq 1$, let $\mathcal{H}_{q}$ be the $q$th Wiener chaos of $X$,
that is, the closed linear subspace of $L^2(\Omega ,\mathcal{F},P)$
generated by the family of random variables $\{H_{q}(X(h)),h\in \EuFrak H,\left\|
h\right\| _{\EuFrak H}=1\}$, where $H_{q}$ is the $q$th Hermite polynomial
defined as
\begin{equation}\label{hermite}
H_q(x) = (-1)^q e^{\frac{x^2}{2}}
 \frac{d^q}{dx^q} \big( e^{-\frac{x^2}{2}} \big).
\end{equation}
The first few Hermite polynomials are $H_1(x)=x$, $H_2(x)=x^2-1$, $H_3(x)=x^3-3x$.
We write by convention $\mathcal{H}_{0} = \mathbb{R}$ and $I_{0}(x)=x$, $x\in\mathbb{R}$. For
any $q\geq 1$, the mapping $I_{q}(h^{\otimes q})=H_{q}(X(h))$ can be extended to a
linear isometry between \ the symmetric tensor product $\EuFrak H^{\odot q}$
equipped with the modified norm
$\left\| \cdot \right\| _{\EuFrak H^{\odot q}}=
\sqrt{q!}
\left\| \cdot \right\| _{\EuFrak H^{\otimes q}}$
and the $q$th Wiener chaos $\mathcal{H}_{q}$.
Then
\begin{equation}\label{e:covcov}
E[I_p(f)I_q(g)]=\delta_{p,q}\times p!\langle f,g\rangle_{\HH^{\otimes p}}
\end{equation}
where $\delta_{p,q}$ stands for the usual Kronecker symbol,
for $f\in\HH^{\odot p}$, $g\in\HH^{\odot q}$ and $p,q\geq 1$.
Moreover, if $f\in\HH^{\otimes q}$, we have
\begin{equation}\label{e:sym}
I_q(f)=I_q(\widetilde{f}),
\end{equation}
where $\widetilde{f}\in\HH^{\odot q}$ is the symmetrization of $f$.

It is well known that $L^2(\Omega,\mathcal{F},P)$ can be decomposed into the infinite
orthogonal sum of the spaces $\mathcal{H}_q$. Therefore, any square integrable random variable
$G\in L^2(\Omega,\mathcal{F},P)$ admits the following Wiener chaotic expansion
\begin{equation}\label{expansion}
G=E[G]+\sum_{q=1}^\infty I_q(f_q),
\end{equation}
where the $f_q\in\HH^{\odot q}$, $q\geq 1$, are uniquely determined by $G$.

Let $\{e_{k},\,k\geq 1\}$ be a complete orthonormal system in $\EuFrak H$.
Given $f\in \EuFrak H^{\odot p}$ and $g\in \EuFrak H^{\odot q}$, for every
$r=0,\ldots ,p\wedge q$, the \textit{contraction} of $f$ and $g$ of order $r$
is the element of $\EuFrak H^{\otimes (p+q-2r)}$ defined by
\begin{equation}
f\otimes _{r}g=\sum_{i_{1},\ldots ,i_{r}=1}^{\infty }\langle
f,e_{i_{1}}\otimes \ldots \otimes e_{i_{r}}\rangle _{\EuFrak H^{\otimes
r}}\otimes \langle g,e_{i_{1}}\otimes \ldots \otimes e_{i_{r}}
\rangle_{\EuFrak H^{\otimes r}}.  \label{v2}
\end{equation}
Since $f\otimes _{r}g$ is not necessarily symmetric, we denote its
symmetrization by $f\widetilde{\otimes }_{r}g\in \EuFrak H^{\odot (p+q-2r)}$.
Observe that $f\otimes _{0}g=f\otimes g$ equals the tensor product of $f$ and
$g$ while, for $p=q$, $f\otimes _{q}g=\langle f,g\rangle _{\EuFrak H^{\otimes q}}$, namely
the scalar product of $f$ and $g$.
In the particular case $\EuFrak H=L^2(A,\mathcal{A},\mu )$, where
$(A,\mathcal{A})$ is a measurable space and $\mu $ is a $\sigma $-finite and
non-atomic measure, one has that $\EuFrak H^{\odot q}=L_{s}^{2}(A^{q},
\mathcal{A}^{\otimes q},\mu ^{\otimes q})$ is the space of symmetric and
square integrable functions on $A^{q}$. In this case, (\ref{v2}) can be rewritten as
\begin{eqnarray*}
(f\otimes _{r}g)(t_1,\ldots,t_{p+q-2r})
&=&\int_{A^{r}}f(t_{1},\ldots ,t_{p-r},s_{1},\ldots ,s_{r}) \\
&&\times \,g(t_{p-r+1},\ldots ,t_{p+q-2r},s_{1},\ldots ,s_{r})d\mu
(s_{1})\ldots d\mu (s_{r}),
\end{eqnarray*}
that is, we identify $r$ variables in $f$ and $g$ and integrate them out.
We shall make use of the following lemma whose proof is a straighforward application
of the definition of contractions and Fubini theorem.

\begin{lemma}\label{contr}
Let $f,g\in\HH^{\odot 2}$. Then
$\|f\otimes_1 g\|^2_{\HH^{\otimes 2}}=
\langle f\otimes_1f,g\otimes_1g\rangle_{\HH^{\otimes 2}}.$
\end{lemma}


\smallskip

Let us now introduce some basic elements of the Malliavin calculus with respect
to the isonormal Gaussian process $X$. Let $\mathcal{S}$
be the set of all
cylindrical random variables of
the form
\begin{equation}
G=\varphi\left( X(h_{1}),\ldots ,X(h_{n})\right) ,  \label{v3}
\end{equation}
where $n\geq 1$, $\varphi:\mathbb{R}^{n}\rightarrow \mathbb{R}$ is an infinitely
differentiable function with compact support and $h_{i}\in \EuFrak H$.
The {\sl Malliavin derivative}  of $G$ with respect to $X$ is the element of
$L^2(\Omega ,\EuFrak H)$ defined as
\begin{equation}\label{e:Dphi}
DG\;=\;\sum_{i=1}^{n}\frac{\partial \varphi}{\partial x_{i}}\left( X(h_{1}),\ldots ,X(h_{n})\right) h_{i}.
\end{equation}
By iteration, one can
define the $m$th derivative $D^{m}G$, which is an element of $L^2(\Omega ,\EuFrak H^{\odot m})$,
for every $m\geq 2$.
For instance, for $G$ as in (\ref{v3}), we have
\[
D^2 G = \sum_{i,j=1}^n \frac{\partial^2 \varphi}{\partial x_i\partial x_j}(X(h_1),\ldots,X(h_n))h_i\otimes h_j.
\]
For $m\geq 1$ and $p\geq 1$, ${\mathbb{D}}^{m,p}$ denotes the closure of
$\mathcal{S}$ with respect to the norm $\Vert \cdot \Vert _{m,p}$, defined by
the relation
\begin{equation}\label{e:norm}
\Vert G\Vert _{m,p}^{p}\;=\;E\left[ |G|^{p}\right] +\sum_{i=1}^{m}E\left(
\Vert D^{i}G\Vert _{\EuFrak H^{\otimes i}}^{p}\right) .
\end{equation}
In particular, $DX(h)=h$ for every $h\in \EuFrak H$. The Malliavin derivative $D$ verifies moreover
the following \textsl{chain rule}. If
$\varphi :\mathbb{R}^{n}\rightarrow \mathbb{R}$ is continuously
differentiable with bounded partial derivatives and if $G=(G_{1},\ldots
,G_{n})$ is a vector of elements of ${\mathbb{D}}^{1,2}$, then $\varphi
(G)\in {\mathbb{D}}^{1,2}$ and
\begin{equation*}
D\varphi (G)=\sum_{i=1}^{n}\frac{\partial \varphi }{\partial x_{i}}(G)DG_{i}.
\end{equation*}

Let now $\HH=
L^{2}(A,\mathcal{A},\mu )$ with $\mu$ non-atomic. Then
an element $u\in\HH$ can be expressed as $u=\{u_t,\,t\in A\}$ and
the Malliavin
derivative of a multiple integral $G$ of the form $I_q(f)$ (with $f\in\HH^{\odot q}$)
is
the element
$DG=\{D_tG,\,t\in A\}$
 of $L^2(A\times \Omega )$ given by
\begin{equation}
D_{t}G=D_t \big[I_q(f)\big]=qI_{q-1}\left( f(\cdot ,t)\right).  \label{dtf}
\end{equation}
Thus the derivative of the random variable $I_q(f)$ is the stochastic process $qI_{q-1}\big(f(\cdot,t)\big)$,
$t\in A$.
Moreover,
\[
\|D\big[I_q(f)\big]\|^2_\HH=q^2\int_A I_{q-1}\left( f(\cdot ,t)\right)^2\mu(dt).
\]
For any $G\in L^2(\Omega,\mathcal{F},P)$
as in (\ref{expansion}), we
define
\begin{equation}\label{e:L-1}
L^{-1}G=-\sum_{q=1}^\infty \frac{1}q I_q(f_q).
\end{equation}
It is proven in \cite{NP07} that for every centered $G\in L^2(\Omega,\mathcal{F},P)$
and every $\mathcal{C}^1$ and Lipschitz function $h:\R\to\C$,
\begin{equation}\label{ipp-np}
E[Gh(G)]=E[h'(G)\langle DG,-DL^{-1}G\rangle_\HH].
\end{equation}
In the particular case $h(x)=x$, we obtain from (\ref{ipp-np}) that
\begin{equation}\label{ipp-nph1}
{\rm Var}[G]=E[G^2]=E[\langle DG,-DL^{-1}G\rangle_\HH],
\end{equation}
where `Var' denotes the variance.
Moreover, if $G\in\mathbb{D}^{2,4}$ is centered, then it is shown in \cite{NPR} that
\begin{equation}\label{variance}
{\rm Var}[\langle DG,-DL^{-1}G\rangle]\leq
\frac52 E[\|DG\|_\HH^4]^{\frac12} E[\|D^2G\otimes_1 D^2G\|_{\HH^{\otimes 2}}^2]^{\frac12}.
\end{equation}
Finally, we shall also use the following bound, established
in a slightly different way in \cite[Corollary 4.2]{NPR},
for the difference between the characteristic functions of a
centered random variable in $\mathbb{D}^{2,4}$ and of a standard Gaussian random variable.

\begin{lemma} \label{noupec}
Let $G\in\mathbb{D}^{2,4}$ be centered.
Then, for any $t\in\R$, we have
\begin{eqnarray}
\big|E[e^{itG}]\!-\!e^{-t^2/2}\big|
\!\leq\! |t|\big|1\!-\!E[G^2]\big|\!+\!\frac{|t|}{2}\sqrt{10}\,
E[\|D^2G\otimes_1D^2G\|^2_{\HH^{\otimes 2}}]^{\frac14} E[\|DG\|^4_\HH]^{\frac14}.
\label{was}
\end{eqnarray}
\end{lemma}
{\bf Proof}. For all $t\in\R$, let $\varphi(t)=e^{t^2/2}
E[e^{itG}]$.
It follows from (\ref{ipp-np}) that
\[
\varphi'(t)=te^{t^2/2}E[e^{itG}]+ie^{t^2/2}E[Ge^{itG}]
=te^{t^2/2}E[e^{itG}(1-\langle DG,-DL^{-1}G\rangle_\HH)].
\]
Hence, we obtain that
\[
\big|\varphi(t)-\varphi(0)\big|\leq
\sup_{u\in[0,\,t]}|\varphi'(u)|\leq |t|e^{t^2/2}E\big[|1-\langle
DG,-DL^{-1}G\rangle_\HH|\big],
\]
which leads to
\[
\big|E[e^{itG}]-e^{-t^2/2}\big|
\leq |t|\,E\big[|1-\langle
DG,-DL^{-1}G\rangle_\HH|\big].
\]
Consequently, we deduce from (\ref{ipp-nph1}) together with Cauchy-Schwarz inequality
that
\begin{eqnarray*}
\big|E[e^{itG}]-e^{-t^2/2}\big|
&\leq&|t|\,\big|1-E[G^2]\big|+|t|\,E\big[|E[G^2]-\langle DG,-DL^{-1}G\rangle_\HH|\big],\\
&\leq&|t|\,\big|1-E[G^2]\big|+|t|\sqrt{{\rm Var}\big(\langle DG,-DL^{-1}G\rangle_\HH\big)}.
\end{eqnarray*}
We conclude the proof of Lemma \ref{noupec} by using (\ref{variance}).
\fin

\section{A criterion for ASCLT on the Wiener space}\label{gen-crit}

The following result, due to Ibragimov and Lifshits \cite{IL},
gives a sufficient condition for extending
convergence in law to ASCLT. It will play a crucial role in all the sequel.

\begin{thm}\label{thm-IL}
Let $\{G_n\}$ be a sequence of random variables converging in distribution towards a random
variable $G_\infty$, and set
\begin{equation}\label{e:delta}
\Delta_n(t)=\frac1{\log n}\sum_{k=1}^n \frac{1}k \big(e^{itG_k}-E(e^{itG_\infty})\big).
\end{equation}
If, for  all $r>0$,
\begin{equation}\label{cond-IL}
\sup_{|t|\leq r}\sum_n \frac{E\vert \Delta_n(t)\vert^2}{n\log n}<\infty,
\end{equation}
then, almost surely, for all continuous and bounded function
$\varphi:\R\to\R$, we have
\[
\frac{1}{\log n}\sum_{k=1}^{n} \frac{1}{k}\,\varphi(G_k)
\longrightarrow E[\varphi(G_\infty)],\quad\mbox{as $n\to\infty$}.
\]
\end{thm}

The following theorem is the main abstract result of this section.
It provides a suitable criterion for an ASCLT for normalized
sequences in $\mathbb{D}^{2,4}$.

\begin{thm}\label{main2}
Let the notation of Section \ref{sec22} prevail.
Let $\{G_n\}$ be a sequence in $\mathbb{D}^{2,4}$ satisfying, for all $n\geq 1$,
$E[G_n]=0$ and $E[G_n^2]=1$.
Assume that
\begin{eqnarray*}
(A_0)&\quad&\displaystyle{\sup_{n\geq 1}E\big[\|DG_n\|^4_\HH]<\infty},
\end{eqnarray*}
and
\[
E[\|D^2G_n\otimes_1D^2G_n\|_{\HH^{\otimes 2}}^2]\to 0,\quad\mbox{as $n\to\infty$}.
\]
Then,
$G_n\overset{\rm law}{\longrightarrow}N\sim
\NN(0,1)$ as $n\to\infty$.
Moreover, assume that the two following conditions also hold
\begin{eqnarray*}
(A_1) &\quad&\displaystyle{\sum_{n\geq 2} \frac1{n\log^2n}\sum_{k=1}^n \frac1k\,
E[\|D^2G_k\otimes_1D^2G_k\|_{\HH^{\otimes 2}}^2]^{\frac14}
<\infty},\\
(A_2) &\quad&\displaystyle{\sum_{n\geq 2} \frac1{n\log^3n}\sum_{k,l=1}^n \frac{\big|
E(G_kG_l)
\big|}
{kl}}<\infty.
\end{eqnarray*}
Then, $\{G_n\}$ satisfies an ASCLT. In other words, almost surely, for all
continuous and bounded function $\varphi:\R\to\R$,
\[
\frac{1}{\log n}\sum_{k=1}^{n} \frac{1}{k}\,\varphi(G_k) \,\longrightarrow\, E[\varphi(N)],\quad\mbox{as $n\to\infty$}.
\]
\end{thm}
\begin{rem}\label{abc-rem}
{\rm
If there exists $\alpha>0$ such that
$E[\|D^2G_k\otimes_1D^2G_k\|_{\HH^{\otimes 2}}^2]=O(k^{-\alpha})$ as $k\to\infty$, then
$(A_1)$ is clearly satisfied.
On the other hand, if there exists $C,\alpha>0$ such that
$\big|E[G_kG_l]\big|\leq C\left(\frac{k}{l}\right)^\alpha$
for all $k\leq l$,
then, for some positive constants $a$, $b$ independent of $n$,
we have
\begin{eqnarray*}
\sum_{n\geq 2}\frac{1}{n\log^3 n}\sum_{l=1}^n\frac{1}{l}\sum_{k=1}^l \frac{\big|E[G_kG_l]\big|}{k}
&\leq&C\sum_{n\geq 2}\frac{1}{n\log^3 n}\sum_{l=1}^n \frac{1}{l^{1+\alpha}}\sum_{k=1}^l k^{\alpha-1},\\
&\leq&
a \sum_{n\geq 2}\frac{1}{n\log^3 n}\sum_{l=1}^n \frac{1}{l}
\leq
b \sum_{n\geq 2}\frac{1}{n\log^2 n}<\infty,
\end{eqnarray*}
which means that $(A_2)$ is also satisfied.
}
\end{rem}
{\bf Proof of Theorem \ref{main2}}.
The fact that $G_n\overset{\rm law}{\longrightarrow}N\sim
\NN(0,1)$ follows from \cite[Corollary 4.2]{NPR}.
In order to prove that the ASCLT holds, we shall verify the sufficient condition (\ref{cond-IL}),
that is the Ibragimov-Lifshits criterion. For simplicity, let $g(t)=E(e^{itN})=e^{-t^2/2}$.
Then, we have
\begin{eqnarray}
\label{train1}
&\!\!&E\vert \Delta_n(t)\vert^2 \\
&\!=\!&\frac{1}{\log^2 n}\sum_{k,l=1}^n \frac{1}{kl}
E\left[\big(e^{itG_k}-g(t)\big)\big(e^{-itG_l}-g(t)\big)\right],\notag\\
&\!=\!&\frac{1}{\log^2 n} \sum_{k,l=1}^{n}\frac{1}{kl}
\left[
E\big(e^{it(G_k-G_l)}\big)-g(t)\left(E\big(e^{itG_k}\big)+E\big(e^{-itG_l}\big)\right)+g^2(t)
\right],\notag\\
&\!=\!&\frac{1}{\log^2 n} \sum_{k,l=1}^{n}\frac{1}{kl}
\left[
\left(E\big(e^{it(G_k-G_l)}\big)-g^2(t)\right)
- g(t)\left(
E\big(e^{itG_k}\big)-g(t)
\right)
- g(t)\left(
E\big(e^{-itG_l}\big)-g(t)
\right)
\right].\notag
\end{eqnarray}
Let $t\in\R$ and $r>0$ be such that $|t|\leq r$.
It follows from inequality (\ref{was}) together with assumption $(A_0)$ that
\begin{equation}
\label{train2}
\left| E\big(e^{itG_k}\big) - g(t)\right| \leq \frac{r\xi}2\,\sqrt{10}\,
E[\|D^2G_k\otimes_1D^2G_k\|_{\HH^{\otimes 2}}^2]^{\frac14}
\end{equation}
where $\xi=\sup_{n\geq 1}E\big[\|DG_n\|^4_\HH]^{\frac14}$. Similarly,
\begin{equation}
\label{train3}
\left| E\big(e^{-itG_l}\big) - g(t)\right| \leq \frac{r\xi}2\,\sqrt{10}\,
E[\|D^2G_l\otimes_1D^2G_l\|_{\HH^{\otimes 2}}^2]^{\frac14}.
\end{equation}
On the other hand, we also have via (\ref{was}) that
\begin{eqnarray*}
&\!\!&\left| E\big(e^{it(G_k-G_l)}\big) - g^2(t)\right|
=\left| E\big(e^{it\sqrt{2}\,\frac{G_k-G_l}{\sqrt{2}}}\big) - g(\sqrt{2}\,t)\right|, \notag\\
&\!\leq\!&  \sqrt{2}r\Bigl|1-\frac12 E[(G_k-G_l)^2]\Bigr|+ r\xi\,\sqrt{5}\,
E[\|D^2(G_k-G_l)\otimes_1 D^2(G_k-G_l)\|^2_{\HH^{\otimes 2}}]^{\frac14},\notag\\
&\!\leq\!&
\sqrt{2}r|E[G_kG_l]|+ r\xi\,\sqrt{5}\,
E[\|D^2(G_k-G_l)\otimes_1 D^2(G_k-G_l)\|^2_{\HH^{\otimes 2}}]^{\frac14}.
\end{eqnarray*}
Moreover
\begin{eqnarray*}
\|D^2(G_k-G_l)\otimes_1 D^2(G_k-G_l)\|^2_{\HH^{\otimes 2}}&\leq&
2\|D^2G_k\otimes_1 D^2G_k\|^2_{\HH^{\otimes 2}}
+
2\|D^2G_l\otimes_1 D^2G_l\|^2_{\HH^{\otimes 2}} \\
&&+ 4\|D^2G_k\otimes_1 D^2G_l\|^2_{\HH^{\otimes 2}}.
\end{eqnarray*}
In addition, we infer from Lemma \ref{contr} that
\begin{eqnarray*}
E\big[
\|D^2G_k\otimes_1 D^2G_l\|^2_{\HH^{\otimes 2}}
\big]&=&E\big[\langle D^2G_k\otimes_1 D^2G_k,D^2G_l\otimes_1D^2G_l\rangle_{\HH^{\otimes 2}}\big],\\
&\leq&\Bigl(E\big[\| D^2G_k\otimes_1 D^2G_k\|^2_{\HH^{\otimes 2}}\big]\Bigr)^{\frac12}
\Bigl(E\big[\| D^2G_l\otimes_1 D^2G_l\|^2_{\HH^{\otimes 2}}\big]\Bigr)^{\frac12},\\
&\leq&\frac12E\big[\| D^2G_k\otimes_1 D^2G_k\|^2_{\HH^{\otimes 2}}\big]+\frac12
E\big[\| D^2G_l\otimes_1 D^2G_l\|^2_{\HH^{\otimes 2}}\big].
\end{eqnarray*}
Consequently, we deduce from the elementary inequality
$(a+b)^{\frac14}\leq a^{\frac14}+b^{\frac14}$ that
\begin{eqnarray}
\label{train4}
&&\left| E\big(e^{it(G_k-G_l)}\big) - g^2(t)\right|\\
&\leq& \sqrt{2} r|E[G_kG_l]|+ r\xi\,\sqrt{10}\,
\Bigl(E\big[
\|D^2G_k\otimes_1 D^2G_k\|^2_{\HH^{\otimes 2}}
\big]^{\frac14}+
E\big[
\|D^2G_l\otimes_1 D^2G_l\|^2_{\HH^{\otimes 2}}
\big]^{\frac14}\Bigr).\notag
\end{eqnarray}
Finally, (\ref{cond-IL}) follows from the conjunction of $(A_1)$ and $(A_2)$ together with
(\ref{train1}), (\ref{train2}), (\ref{train3}) and
(\ref{train4}), which completes the proof of Theorem \ref{main2}.
\fin

We now provide an explicit application of Theorem \ref{main2}.
\begin{thm}\label{ex-thm}
Let $X=\{X_n\}_{n\in\mathbb{Z}}$ denote a centered stationary Gaussian
sequence with unit variance, such that $\sum_{r\in\mathbb{Z}}|\rho(r)|<\infty$, where
$\rho(r)=E[X_0X_{r}]$.
Let $f:\R\to\R$ be a symmetric real function of class $\mathcal{C}^2$,
and let $N\sim \mathscr{N}(0,1)$.
Assume moreover that $f$ is not constant and that $E[f''(N)^4]<\infty$.
For any $n\geq 1$, let
\[
G_n=\frac{1}{\sigma_n\sqrt{n}}\sum_{k=1}^{n} \big(f(X_k)-E[f(X_k)]\big)
\]
where $\sigma_n$ is the positive normalizing constant which ensures that $E[G_n^2]=1$.
Then, as $n\to\infty$,
$G_n\overset{\rm law}{\longrightarrow} N$
and $\{G_n\}$ satisfies an ASCLT.
In other words, almost surely, for all
continuous and bounded function $\varphi:\R\to\R$,
\[
\frac{1}{\log n}\sum_{k=1}^{n} \frac{1}{k}\,\varphi(G_k) \,\longrightarrow\, E[\varphi( N)],\quad\mbox{as $n\to\infty$}.
\]
\end{thm}
\begin{rem}\label{tyty}
{\rm
We can replace the assumption `$f$ is symmetric and non-constant'
by
\[
\sum_{q=1}^\infty \frac{1}{q!}\big(E[f(N)H_q(N)]\big)^2\sum_{r\in\mathbb{Z}} |\rho(r)|^{q}<\infty\,\mbox{ and }\,
\sum_{q=1}^\infty \frac{1}{q!}\big(E[f(N)H_q(N)]\big)^2\sum_{r\in\mathbb{Z}} \rho(r)^{q}>0.
\]
Indeed, it suffices to replace the monotone convergence argument used to prove (\ref{bgbc}) below
by a bounded convergence argument.
However, this new assumption seems rather difficult to check in general, except of course
when the sum with respect to $q$ is finite, that is when $f$ is a polynomial.
}
\end{rem}
{\bf Proof of Theorem \ref{ex-thm}}. First, note that a consequence
of \cite[inequality (3.19)]{NPR} is that we automatically
have $E[f'(N)^4]<\infty$ and $E[f(N)^4]<\infty$.
Let us now expand $f$ in terms of Hermite polynomials. Since $f$ is symmetric, we can write
\[
f=E[f(N)]+\sum_{q=1}^\infty c_{2q}H_{2q},
\]
where the real numbers $c_{2q}$ are given by $(2q)!c_{2q}=E[f(N)H_{2q}(N)]$. Consequently,
\begin{eqnarray*}
\sigma_n^2&=&\frac{1}{n}\sum_{k,l=1}^n {\rm Cov}[f(X_k),f(X_l)]
=\sum_{q=1}^\infty c_{2q}^2(2q)!\,\frac{1}{n}\sum_{k,l=1}^n \rho(k-l)^{2q},\\
&=&\sum_{q=1}^\infty c_{2q}^2(2q)!\sum_{r\in\mathbb{Z}} \rho(r)^{2q}
\left(1-\frac{|r|}{n}\right){\bf 1}_{\{
|r|\leq n\}}.
\end{eqnarray*}
Hence, it follows from the monotone convergence theorem that
\begin{equation}
\sigma_n^2\,\longrightarrow\, \sigma_\infty^2=\sum_{q=1}^\infty c_{2q}^2(2q)!\sum_{r\in\mathbb{Z}} \rho(r)^{2q},\quad\mbox{as $n\to\infty$}.
\label{bgbc}
\end{equation}
Since $f$ is not constant, one can find some $q\geq 1$ such that
$c_{2q}\neq 0$. Moreover, we also have $\sum_{r\in\mathbb{Z}}\rho(r)^{2q}\geq \rho(0)^{2q}=1$.
Hence, $\sigma_\infty> 0$, which implies in particular that the infimum of the sequence
$\{\sigma_n\}_{n\geq 1}$ is positive.

The Gaussian space generated by $X=\{X_k\}_{k\in\mathbb{Z}}$ can be identified with an isonormal Gaussian
process of the type $X=\{X(h):\,h\in\HH\}$, for $\HH$ defined as follows:
(i) denote by $\mathcal{E}$ the set of all sequences indexed by $\mathbb{Z}$ with finite support;
(ii) define $\HH$ as the Hilbert space obtained by closing $\mathcal{E}$ with respect to the
scalar product
\begin{equation}\label{e:sproduct}
\langle u,v\rangle_\HH = \sum_{k,l\in\mathbb{Z}} u_kv_l \rho(k-l).
\end{equation}
In this setting, we have $X(\varepsilon_k)=X_k$ where $\varepsilon_k=\{\delta_{kl}\}_{l\in\mathbb{Z}}$, $\delta_{kl}$ standing for the Kronecker symbol.
In view of (\ref{e:Dphi}), we have
\[
DG_n=\frac{1}{\sigma_n\sqrt{n}}\sum_{k=1}^n f'(X_k)\varepsilon_k.
\]
Hence
\[
\|DG_n\|^2_\HH = \frac{1}{\sigma_n^2\,n}\sum_{k,l=1}^n f'(X_k)f'(X_l)\langle \e_k,\e_l\rangle_\HH
=\frac{1}{\sigma_n^2\,n}\sum_{k,l=1}^n f'(X_k)f'(X_l)\rho(k-l),
\]
so that
\[
\|DG_n\|^4_\HH =\frac{1}{\sigma_n^4\,n^2}\sum_{i,j,k,l=1}^n
f'(X_i)f'(X_j)f'(X_k)f'(X_l)
\rho(i-j)\rho(k-l).
\]
We deduce from Cauchy-Schwarz inequality that
\[
\big|E[f'(X_i)f'(X_j)f'(X_k)f'(X_l)]\big|\leq (E[f'(N)^4])^\frac14,
\]
which leads to
\begin{equation}\label{co1}
E[\|DG_n\|^4_\HH]\leq \frac{1}{\sigma_n^4}\Bigl(E[f'(N)^4]\Bigr)^{\frac14}\left(\sum_{r\in\mathbb{Z}}|\rho(r)|\right)^2.
\end{equation}
On the other hand, we also have
\[
D^2G_n = \frac{1}{\sigma_n\sqrt{n}}\sum_{k=1}^n f''(X_k)\varepsilon_k\otimes \varepsilon_k,
\]
and therefore
\[
D^2G_n\otimes_1 D^2G_n = \frac{1}{\sigma_n^2\,n}\sum_{k,l=1}^n f''(X_k)f''(X_l)\rho(k-l)\varepsilon_k\otimes
\varepsilon_l.
\]
Hence
\begin{eqnarray}
&&E\big[\|D^2G_n\otimes_1 D^2G_n\|^2_{\HH^{\otimes 2}}\big],\notag\\
& =& \frac{1}{\sigma_n^4\,n^2}
\sum_{i,jk,l=1}^n E\big[f''(X_i)f''(X_j)f''(X_k)f''(X_l)\big]\rho(k-l)\rho(i-j)\rho(k-i)\rho(l-j),\notag\\
&\leq &\frac{(E\big[f''(N)^4\big])^\frac14}{\sigma_n^4\,n}\sum_{u,v,w\in\mathbb{Z}}
|\rho(u)||\rho(v)||\rho(w)||\rho(-u+v+w)|,\notag\\
&\leq &\frac{(E\big[f''(N)^4\big])^\frac14\|\rho\|_\infty}{\sigma_n^4\,n}\left(\sum_{r\in\mathbb{Z}}
|\rho(r)|\right)^3<\infty.\label{co2}
\end{eqnarray}
By virtue of Theorem \ref{main2} together with the fact that
$\inf_{n\geq 1}\sigma_n>0$, the inequalities (\ref{co1}) and (\ref{co2}) imply
that $G_n\overset{\rm law}{\longrightarrow} N$.
Now, in order to show that the ASCLT holds, we shall also check that
conditions $(A_1)$ and
$(A_2)$ in Theorem \ref{main2} are fulfilled. First, still because
$\inf_{n\geq 1}\sigma_n>0$, $(A_1)$ holds
since we have
$E\big[\|D^2G_n\otimes_1 D^2G_n\|^2_{\HH^{\otimes 2}}\big]=O(n^{-1})$
by (\ref{co2}), see also Remark \ref{abc-rem}. Therefore, it only remains
to prove $(A_2)$. Gebelein's inequality (see e.g. identity (1.7) in \cite{gebelein}) states that
\[
\big|{\rm Cov}[f(X_i),f(X_j)]\big|\leq E[X_iX_j]\sqrt{{\rm Var}[f(X_i)]}\sqrt{{\rm Var}[f(X_j)]}
=\rho(i-j){\rm Var}[f(N)].
\]
Consequently,
\begin{eqnarray*}
 \big|E[G_kG_l]\big|&=&\frac{1}{\sigma_k\sigma_l\sqrt{kl}}\left|
\sum_{i=1}^k\sum_{j=1}^l {\rm Cov}[f(X_i),f(X_j)]\right|
\leq
\frac{{\rm Var}[f(N)]}{\sigma_k\sigma_l\sqrt{kl}}
\sum_{i=1}^k\sum_{j=1}^l |\rho(i-j)|,\\
&=&
\frac{{\rm Var}[f(N)]}{\sigma_k\sigma_l\sqrt{kl}}
\sum_{i=1}^{k}\sum_{r=i-l}^{i-1}|\rho(r)|
\leq
\frac{{\rm Var}[f(N)]}{\sigma_k\sigma_l}\,\sqrt{\frac{k}l}\sum_{r\in\mathbb{Z}}|\rho(r)|.
\end{eqnarray*}
Finally, via the same arguments as in Remark \ref{abc-rem}, $(A_2)$ is satisfied, which completes
the proof of Theorem \ref{ex-thm}.
\fin

The following result specializes Theorem \ref{main2}, by providing a criterion for an ASCLT for
{\sl multiple stochastic integrals} of fixed order $q\geq 2$.
It is expressed in terms of the kernels of these integrals.
\begin{cor}\label{main}
Let the notation of Section \ref{sec22} prevail.
Fix $q\geq 2$, and let $\{G_n\}$ be a sequence of the form $G_n=I_q(f_n)$, with $f_n\in\EuFrak H^{\odot q}$. Assume that
$E[G_n^2]=q!\|f_n\|^2_{\HH^{\otimes q}}=1$ for all $n$, and that
\begin{equation}\label{e:NP}
\|f_n\otimes_r f_n\|_{\EuFrak H^{\otimes 2(q-r)}}\to 0\quad\mbox{ as $n\to\infty$,}\quad\mbox{for every $r=1,\ldots,q-1$.}
\end{equation}
Then,
$G_n\overset{\rm law}{\longrightarrow}N\sim
\NN(0,1)$ as $n\to\infty$.
Moreover, if the two following conditions are also satisfied
\[
(A'_1) \quad\quad\displaystyle{\sum_{n\geq 2} \frac1{n\log^2n}\sum_{k=1}^n \frac1k\,
\|f_k\otimes_r f_k\|_{\EuFrak H^{\otimes 2(q-r)}}}
<\infty\quad\mbox{for every $r=1,\ldots,q-1$},
\]
\[
(A'_2) \quad\quad\displaystyle{\sum_{n\geq 2} \frac1{n\log^3n}\sum_{k,l=1}^n \frac{\big|
\langle f_k,f_l\rangle_{\EuFrak H^{\otimes q}}
\big|}
{kl}}<\infty.
\]
then $\{G_n\}$ satisfies an ASCLT. In other words, almost surely, for all
continuous and bounded function $\varphi:\R\to\R$,
\[
\frac{1}{\log n}\sum_{k=1}^{n} \frac{1}{k}\,\varphi(G_k) \,
\longrightarrow\, E[\varphi(N)],\quad\mbox{as $n\to\infty$}.
\]
\end{cor}
{\bf Proof of Corollary \ref{main}}.
The fact that
$G_n\overset{\rm law}{\longrightarrow}N\sim\NN(0,1)$
follows directly from (\ref{e:NP}), which is the Nualart-Peccati \cite{NP} criterion of normality. In order to prove that
the ASCLT holds, we shall apply once again Theorem \ref{main2}. This is possible because a multiple
integral is always an element of $\mathbb{D}^{2,4}$.
We have, by (\ref{ipp-nph1}),
\[
1=E[G_k^2]=E[\langle DG_k,-DL^{-1}G_k\rangle_\HH]
=\frac{1}q\,E[\|DG_k\|^2_\HH],
\]
where the last inequality follows from $-L^{-1}G_k=\frac1q G_k$, using the definition (\ref{e:L-1}) of $L^{-1}$.
In addition, as the random variables $\|DG_k\|^2_\HH$
live inside the finite sum of the  first $2q$ Wiener chaoses (where all the $L^p$ norm
are equivalent), we deduce that condition $(A_0)$ of Theorem \ref{main2} is satisfied.
On the other hand, it is proven in \cite[page 604]{NPR} that
\[
E\big[\|D^2G_k\otimes_1 D^2G_k\|^2_{\HH^{\otimes 2}}\big]
\leq q^4(q-1)^4\sum_{r=1}^{q-1}(r-1)!^2\binom{q-2}{r-1}^4(2q-2-2r)!\|f_k\otimes_rf_k\|^2_{\HH^{\otimes 2(q-r)}}.
\]
Consequently, condition $(A'_1)$ implies condition $(A_1)$ of Theorem \ref{main2}.
Furthermore, by (\ref{e:covcov}),
$E[G_kG_l]=E\big[I_q(f_k)I_q(f_l)\big]=q!\langle f_k,f_l\rangle_{\EuFrak H^{\otimes q}}$.
Thus, condition $(A'_2)$ is equivalent to condition $(A_2)$ of Theorem \ref{main2},
and the proof of the corollary is done.
\fin

In Corollary \ref{main}, we supposed $q\geq 2$, which implies that $G_n=I_q(f_n)$ is a multiple integral of order at least 2 and hence is not Gaussian.
We now consider the Gaussian case $q=1$.
\begin{cor}\label{q=1}
Let $\{G_n\}$ be a centered Gaussian sequence with unit variance. If the
condition $(A_2)$ in Theorem \ref{main2} is
satisfied,
then $\{G_n\}$ satisfies an ASCLT. In other words, almost surely, for all
continuous and bounded function $\varphi:\R\to\R$,
\[
\frac{1}{\log n}\sum_{k=1}^{n} \frac{1}{k}\,\varphi(G_k) \,
\longrightarrow\, E[\varphi(N)],\quad\mbox{as $n\to\infty$}.
\]
\end{cor}
{\bf Proof of Corollary \ref{q=1}}. Let $t\in\R$ and $r>0$ be such that $|t|\leq r$, and let
$\Delta_n(t)$ be defined as in (\ref{e:delta}).
We have
\begin{eqnarray*}
E\vert \Delta_n(t)\vert^2 &=&\frac{1}{\log^2 n}\sum_{k,l=1}^n \frac{1}{kl}
E\left[\big(e^{itG_k}-e^{-t^2/2}\big)\big(e^{-itG_l}-e^{-t^2/2}\big)\right],\\
&=&\frac{1}{\log^2 n} \sum_{k,l=1}^{n}\frac{1}{kl}
\left[
E\big(e^{it(G_k-G_l)}\big)-e^{-t^2}
\right],\\
&=&\frac{1}{\log^2 n} \sum_{k,l=1}^{n}\frac{e^{-t^2}}{kl}
\big(e^{E(G_kG_l)t^2}-1\big),\\
&\leq&\frac{r^2e^{r^2}}{\log^2 n} \sum_{k,l=1}^{n}\frac{\big|E(G_kG_l)\big|}{kl},
\end{eqnarray*}
since $|e^{x}-1|\leq e^{|x|}|x|$ and $|E(G_kG_l)|\leq 1$.
Therefore, assumption $(A_2)$ implies (\ref{cond-IL}), and the proof of the corollary is done.
\fin

\section{Application to discrete-time fractional Brownian motion}\label{s:FBM}

Let us apply Corollary \ref{q=1} to the particular case $G_n=B^H_n/n^H$,
where $B^H$ is a fractional Brownian motion
with Hurst index $H\in(0,1)$.
We recall that $B^H=(B^H_t)_{t\geq 0}$ is a centered Gaussian process with continuous paths such that
\[
E[B^H_tB^H_s]=\frac{1}{2}\Bigl(t^{2H}+s^{2H} -|t-s|^{2H}\Bigr),\quad s,t\geq 0.
\]
The process $B^H$ is self-similar with stationary increments
and we refer the reader to Nualart \cite{Nbook} and Samorodnitsky and Taqqu \cite{SamorTaqqu}
for its main properties. The increments
\[
Y_k=B^H_{k+1}-B^H_{k},\quad k\geq 0,
\]
called `fractional Gaussian noise',
are centered stationary Gaussian random variables with covariance
\begin{equation}\label{e:cov}
\rho(r)=E[Y_kY_{k+r}]=\frac12\big(|r+1|^{2H}+|r-1|^{2H}-2|r|^{2H}\big), \quad r\in\mathbb{Z}.
\end{equation}
This covariance behaves asymptotically as
\begin{equation}
\label{behaviorrho}
\rho(r)\sim H(2H-1)|r|^{2H-2}\quad\mbox{as $|r|\to\infty$}.
\end{equation}
Observe that $\rho(0)=1$ and
\begin{itemize}
\item[1)] For $0<H<1/2$, $\rho(r)<0$ for $r\neq 0$,
$$\sum_{r\in\Z}|\rho(r)|<\infty
\hspace{1cm}\text{and}\hspace{1cm}
\sum_{r\in\Z}\rho(r)=0.$$
\item[2)] For $H=1/2$, $\rho(r)=0$ if $r\neq 0$.
\item[3)] For $1/2<H<1$,
$$\sum_{r\in\Z}|\rho(r)|=\infty.$$
\end{itemize}
The Hurst index measures the strenght of the dependence when $H\geq 1/2$: the larger $H$ is,
the stronger is the dependence.

A continuous
time version of the following result was obtained by
Berkes and Horv\'ath \cite{BerkesHorvath}
via a different approach.

\begin{thm}\label{thm1234}
For all $H\in(0,1)$, we have,
almost surely, for all continuous and bounded function
$\varphi:\R\to\R$,
\[
\frac1{\log n}\sum_{k=1}^{n}\frac1k\,\varphi(B^H_k/k^H)\longrightarrow E[\varphi(N)],\quad\mbox{as $n\to\infty$}.
\]
\end{thm}
{\bf Proof of Theorem \ref{thm1234}}.
We shall make use of Corollary \ref{q=1}. The cases $H<1/2$ and $H\geq 1/2$ are treated separately.
From now on, the value of a constant $C>0$ may change from line to
line, and we set $\rho(r)=\frac12\big(|r+1|^{2H}+|r-1|^{2H}-2|r|^{2H}\big)$, $r\in\mathbb{Z}$.\\

{\it Case $H<1/2$.} For any $b\geq a\geq 0$, we have
\[
b^{2H}-a^{2H}=2H \int_0^{b-a}\frac{dx}{(x+a)^{1-2H}}
\leq 2H\int_0^{b-a}\frac{dx}{x^{1-2H}}=(b-a)^{2H}.
\]
Hence, for $l\geq k\geq 1$, we have $l^{2H}-(l-k)^{2H}\leq k^{2H}$ so that
\[
|E[B^H_kB^H_l]|=\frac12\big(k^{2H}+l^{2H}-(l-k)^{2H}\big)\leq k^{2H}.
\]
Thus
\begin{eqnarray*}
\sum_{n\geq 2}\frac1{n\log^3 n}\sum_{l=1}^n\frac1l\sum_{k=1}^l\frac{|E[G_kG_l]|}{k}
&\!=\!&\sum_{n\geq 2}\frac1{n\log^3 n}\sum_{l=1}^n\frac1{l^{1+H}}\sum_{k=1}^l\frac{|E[B^H_kB^H_l]|}{k^{1+H}},\\
&\!\leq\!&\sum_{n\geq 2}\frac1{n\log^3 n}\sum_{l=1}^n\frac1{l^{1+H}}\sum_{k=1}^l\frac1{k^{1-H}},\\
&\!\leq\!&C\sum_{n\geq 2}\frac1{n\log^3 n}\sum_{l=1}^n\frac1{l}
\leq C\sum_{n\geq 2}\frac1{n\log^2 n}<\infty.
\end{eqnarray*}
Consequently, condition $(A_2)$ in Theorem \ref{main2} is satisfied.\\

{\it Case $H\geq 1/2$}. For $l\geq k\geq 1$, it follows from (\ref{e:cov})-(\ref{behaviorrho})
that
\begin{eqnarray*}
|E[B^H_kB^H_l]|&=&\left|
\sum_{i=0}^{k-1}\sum_{j=0}^{l-1}
E[(B^H_{i+1}-B^H_i)(B^H_{j+1}-B^H_j)]
\right|
\leq\sum_{i=0}^{k-1}\sum_{j=0}^{l-1}|\rho(i-j)|,\\
&\leq&k\sum_{r=-l+1}^{l-1}|\rho(r)|\leq Ckl^{2H-1}.
\end{eqnarray*}
The last inequality comes from the fact that
$\rho(0)=1$, $\rho(1)=\rho(-1)=(2^{2H}-1)/2$
and, if $r\geq 2$,
\begin{eqnarray*}
|\rho(-r)|&=&|\rho(r)|=\big| E[(B^H_{r+1}-B^H_r)B^H_1]
=H(2H-1)\int_0^1 du\int_r^{r+1}dv(v-u)^{2H-2}\\
&\leq& H(2H-1)\int_0^1 (r-u)^{2H-2}du \leq H(2H-1)(r-1)^{2H-2}.
\end{eqnarray*}
Consequently,
\begin{eqnarray*}
\sum_{n\geq 2}\frac1{n\log^3 n}\sum_{l=1}^n\frac1l\sum_{k=1}^l\frac{|E[G_kG_l]|}{k}
&=&\sum_{n\geq 2}\frac1{n\log^3 n}\sum_{l=1}^n\frac1{l^{1+H}}\sum_{k=1}^l\frac{|E[B^H_kB^H_l]|}{k^{1+H}},\\
&\leq&C\sum_{n\geq 2}\frac1{n\log^3 n}\sum_{l=1}^n\frac1{l^{2-H}}\sum_{k=1}^l\frac1{k^H},\\
&\leq&C\sum_{n\geq 2}\frac1{n\log^3 n}\sum_{l=1}^n\frac1{l}
\leq C\sum_{n\geq 2}\frac1{n\log^2 n}<\infty.
\end{eqnarray*}
Finally, condition $(A_2)$ in Theorem \ref{main2} is satisfied,
which completes the proof of Theorem \ref{thm1234}.
\fin

\section{Partial sums of Hermite polynomials: the Gaussian limit case}

Let $X=\{X_k\}_{k\in\mathbb{Z}}$ be a centered stationary Gaussian
process and for all $r\in\mathbb{Z}$, set $\rho(r)=E[X_0X_{r}]$.
Fix an integer $q\geq 2$, and let $H_q$ stands for the Hermite polynomial of degree $q$, see
(\ref{hermite}).
We are interested in an ASCLT for the $q$-Hermite power variations of $X$,  defined as
\begin{equation}\label{her-intro1}
V_n=\sum_{k=1}^{n} H_q(X_k),\quad n\geq 1,
\end{equation}
in cases where $V_n$, adequably normalized, converges to a normal distribution.
Our result is as follows.
\begin{thm}\label{BM}
Assume that $\sum_{r\in\mathbb{Z}}|\rho(r)|^q<\infty$, that
$\sum_{r\in\mathbb{Z}}\rho(r)^q>0$ and that there exists $\alpha>0$ such
that $\sum_{|r|>n}|\rho(r)|^q=O(n^{-\alpha})$, as $n\to\infty$.
For any $n\geq 1$, define
\[
G_n=\frac{V_n}{\sigma_n\,\sqrt{n}},
\]
where $V_n$ is given by (\ref{her-intro1}) and
$\sigma_n$ denotes the positive normalizing constant
which ensures that $E[G_n^2]=1$.
Then
$G_n\overset{\rm law}{\longrightarrow}N\sim
\NN(0,1)$ as $n\to\infty$, and
$\{G_n\}$ satisfies an ASCLT. In other words,
almost surely, for all continuous and bounded function $\varphi:\R\to\R$,
\[
\frac{1}{\log n}\sum_{k=1}^{n} \frac{1}{k}\,\varphi(G_k)
\longrightarrow E[\varphi(N)],\quad\mbox{as $n\to\infty$}.
\]
\end{thm}
{\it Proof}. We shall make use of Corollary \ref{main}. Let $C$ be a positive
constant, depending only on $q$ and $\rho$, whose value may change from line to line.
We consider the real and separable Hilbert space $\EuFrak H$ as defined in the
proof of Theorem \ref{ex-thm}, with the scalar product (\ref{e:sproduct}). Following the same line of reasoning as in
the proof of (\ref{bgbc}),
it is possible to show that
$\sigma_n^2\to q!\sum_{r\in\mathbb{Z}}\rho(r)^q>0$. In particular,
the infimum of the sequence
$\{\sigma_n\}_{n\geq 1}$ is positive. On the other hand, we have $G_n=I_q(f_n)$,
where the kernel $f_n$ is given by
\[
f_n=\frac{1}{\sigma_n\sqrt{n}}\sum_{k=1}^n \varepsilon_k^{\otimes q},
\]
with $\varepsilon_k=\{\delta_{kl}\}_{l\in\mathbb{Z}}$, $\delta_{kl}$ standing for the Kronecker symbol.
For all $n\geq 1$ and $r=1,\ldots,q-1$, we have
\[
f_n\otimes_r f_n = \frac{1}{\sigma_n^2\,n}
\sum_{k,l=1}^n\rho(k-l)^r\varepsilon_k^{\otimes(q-r)}\otimes\varepsilon_l^{\otimes(q-r)}.
\]
We deduce that
\[
\|f_n\otimes_r f_n\|^2_{\HH^{\otimes(2q-2r)}}=\frac{1}{\sigma_n^4n^2}
\sum_{i,j,k,l=1}^n \rho(k-l)^r\rho(i-j)^r\rho(k-i)^{q-r}\rho(l-j)^{q-r}.
\]
Consequently, as in the proof of (\ref{co2}), we obtain that
$\|f_n\otimes_r f_n\|^2_{\HH^{\otimes(2q-2r)}}\leq A_n$ where
\[
A_n=\frac{1}{\sigma_n^4n}
\sum_{u,v,w\in D_n} |\rho(u)|^r|\rho(v)|^r|\rho(w)|^{q-r}|\rho(-u+v+w)|^{q-r}
\]
with $D_n=\{-n,\ldots,n\}$. Fix an integer $m\geq 1$ such that $n\geq m$.
We can split $A_n$ into two terms $A_n=B_{n,m}+C_{n,m}$ where
\begin{eqnarray*}
B_{n,m}&\!=\!&\frac{1}{\sigma_n^4n}
\sum_{u,v,w\in D_m}
|\rho(u)|^r|\rho(v)|^r|\rho(w)|^{q-r}|\rho(-u+v+w)|^{q-r},\\
C_{n,m}&\!=\!&\frac{1}{\sigma_n^4n}
\sum_{\substack{u,v,w\in D_n\\|u|\vee|v|\vee|w|> m}}
|\rho(u)|^r|\rho(v)|^r|\rho(w)|^{q-r}|\rho(-u+v+w)|^{q-r}.\\
\end{eqnarray*}
We clearly have
\[
B_{n,m}\leq \frac{1}{\sigma_n^4n}\|\rho\|_\infty^{2q}(2m+1)^3\leq \frac{Cm^3}{n}.
\]
On the other hand,
$D_n\cap\{|u|\vee|v|\vee|w|>m\}\subset D_{n,m,u}\cup D_{n,m,v}\cup D_{n,m,w}$
where the set $D_{n,m,u}=\{|u|>m,|v|\leq n,|w|\leq n\}$ and a similar definition for
$D_{n,m,v}$ and $D_{n,m,w}$. Denote
\[
C_{n,m,u}=\frac{1}{\sigma_n^4n}
\sum_{u,v,w\in D_{n,m,u}}
 |\rho(u)|^r|\rho(v)|^r|\rho(w)|^{q-r}|\rho(-u+v+w)|^{q-r}
\]
and a similar expression for
$C_{n,m,v}$ and $C_{n,m,w}$.
It follows from H\"older inequality that
\begin{equation}
\label{holdercnmu}
C_{n,m,u}\leq \frac{1}{\sigma_n^4n}
\!\left(\!\sum_{u,v,w\in D_{n,m,u}}\!\!\!\!|\rho(u)|^q|\rho(v)|^q \!\right)^{\frac{r}{q}}\!\!
\left(\!\sum_{u,v,w\in D_{n,m,u}}\!\!\!\!|\rho(w)|^q|\rho(-u+v+w)|^q\!\right)^{1-\frac{r}{q}}\!\!\!\!\!.
\end{equation}
However,
\[
\sum_{u,v,w\in D_{n,m,u}}|\rho(u)|^q|\rho(v)|^q\leq (2n+1)
\sum_{|u|>m}|\rho(u)|^q \sum_{v\in\mathbb{Z}}|\rho(v)|^q
\leq Cn \sum_{|u|>m}|\rho(u)|^q.
\]
Similarly,
\[
\sum_{u,v,w\in D_{n,m,u}}|\rho(w)|^q|\rho(-u+v+w)|^q\leq (2n+1)
\sum_{v\in\mathbb{Z}}|\rho(v)|^q \sum_{w\in\mathbb{Z}}|\rho(w)|^q
\leq Cn.
\]
Therefore, (\ref{holdercnmu}) and the last assumption of Theorem \ref{BM}
imply that for $m$ large enough
\begin{equation*}
C_{n,m,u}\leq C
\left(\sum_{|u|>m}|\rho(u)|^q\right)^{\frac{r}{q}}\leq Cm^{-\frac{\alpha r}q}.
\end{equation*}
We obtain exactly the same bound for $C_{n,m,v}$ and $C_{n,m,w}$. Combining all these estimates, we finally
find that
\[
\|f_n\otimes_r f_n\|^2_{\HH^{\otimes(2q-2r)}}\leq C \times \inf_{m\leq n}
\left\{\frac{m^3}{n}+
m^{-\frac{\alpha r}q}
\right\}\leq Cn^{-\frac{\alpha r}{3q+\alpha r}}
\]
by taking the value $m=n^{\frac{q}{3q+\alpha r}}$. It
ensures that condition $(A'_1)$ in Corollary \ref{main} is met.
Let us now prove $(A'_2)$. We have
\begin{eqnarray*}
\langle f_k,f_l\rangle_{\EuFrak H^{\otimes q}}
&=&\frac{1}{\sigma_k\sigma_l\sqrt{kl}}\left|
\sum_{i=1}^k\sum_{j=1}^l \rho(i-j)^q \right|
\leq
\frac{1}{\sigma_k\sigma_l\sqrt{kl}}
\sum_{i=1}^k\sum_{j=1}^l |\rho(i-j)|^q,\\
&\leq&
\frac{1}{\sigma_k\sigma_l}\,\sqrt{\frac{k}l}\sum_{r\in\mathbb{Z}}|\rho(r)|^q,
\end{eqnarray*}
so $(A'_2)$ is also satisfied, see Remark \ref{abc-rem},
which completes the proof of Theorem \ref{BM}. \fin
The following result contains an explicit situation where the assumptions in Theorem \ref{BM} are in order.
\begin{prop}\label{prop-az}
Assume that
$\rho(r)\sim |r|^{-\beta}L(r)$, as $|r|\to\infty$, for some $\beta>1/q$ and some slowly varying function $L$.
Then $\sum_{r\in\mathbb{Z}}|\rho(r)|^q<\infty$ and there exists $\alpha>0$ such
that $\sum_{|r|>n}|\rho(r)|^q=O(n^{-\alpha})$, as $n\to\infty$.
\end{prop}
{\it Proof.}
By a Riemann sum argument, it is immediate that
$\sum_{r\in\mathbb{Z}}|\rho(r)|^q<\infty$.
Moreover, by \cite[Prop. 1.5.10]{BGT}, we have $\sum_{|r|>n}|\rho(r)|^q\sim
\frac{2}{\beta q-1}n^{1-\beta q}L^q(n)$ so that we can choose $\alpha=\frac12(\beta q-1)>0$ (for instance).
\fin

\section{Partial sums of Hermite polynomials of increments of fractional Brownian motion}\label{s:hermite}

We focus here on increments of the fractional Brownian motion $B^H$ (see Section 4 for details about $B^H$).
More precisely, for every $q\geq 1$, we are interested in an ASCLT for
the $q$-Hermite power variation of $B^H$, defined as
\begin{equation}\label{her-intro}
V_n=\sum_{k=0}^{n-1} H_q(B^H_{k+1}-B^H_{k}),\quad n\geq 1,
\end{equation}
where $H_q$ stands for the Hermite polynomial of degree $q$ given by
(\ref{hermite}).
Observe that Theorem \ref{thm1234} corresponds to the particular case $q=1$. That is why,
from now on, we assume that $q\geq 2$.
When $H\neq 1/2$, the increments of $B^H$ are not independent, so
the asymptotic behavior of \eqref{her-intro} is difficult to investigate because $V_n$ is not linear.
In fact, thanks to the seminal works of
Breuer and Major \cite{BM},
Dobrushin and Major \cite{DoMa},
Giraitis and Surgailis \cite{GS}
and Taqqu \cite{T}, it is  known (recall that $q\geq 2$) that,
as $n\to\infty$
\begin{itemize}
\item If $0<H<1-\frac{1}{2q}$, then
\begin{equation}\label{eq:Breuer_Major1}
G_n:=\frac{V_n}{ \sigma_n\,\sqrt n}
\,\overset{{\rm law}}{\longrightarrow}\,
\NN(0,1).
\end{equation}
\item If $H=1-\frac{1}{2q}$, then
\begin{equation}
\label{eq:Breuer_Major2}
G_n:=\frac{V_n}{\sigma_n\sqrt{ n \log n}}
\,\overset{{\rm law}}{\longrightarrow}\,
\NN(0,1).
\end{equation}
\item If $H>1-\frac{1}{2q}$, then
\begin{equation}
\label{eq:Breuer_Major3}
G_n:=n^{q(1-H)-1}V_n
\,\overset{{\rm law}}{\longrightarrow}G_\infty
\end{equation}
\end{itemize}
where $G_\infty$ has an `Hermite distribution'.
Here, $\sigma_n$ denotes the positive normalizing constant
which ensures that $E[G_n^2]=1$. The proofs of (\ref{eq:Breuer_Major1}) and (\ref{eq:Breuer_Major2}), together with rates of convergence,
can be found in \cite{NP07} and \cite{BN}, respectively.
A short proof of (\ref{eq:Breuer_Major3}) is given in Proposition \ref{lm-hermite-yahoo} below.
Notice that rates of convergence can be found in \cite{BN}.
Our proof of (\ref{eq:Breuer_Major3})  is based on the fact that, for {\it fixed} $n$,
$Z_n$ defined in (\ref{sn}) below and $G_n$
share the same law,
because of the self-similarity property of fractional Brownian motion.
\begin{prop}\label{lm-hermite-yahoo}
Assume $H>1-\frac{1}{2q}$, and define $Z_n$
by
\begin{equation}\label{sn}
Z_n=n^{q(1-H)-1}\sum_{k=0}^{n-1} H_q\big(n^H(B^H_{(k+1)/n}-B^H_{k/n})\big),\quad n\geq 1.
\end{equation}
Then, as $n\to\infty$,
$\{Z_n\}$ converges almost surely and in $L^2(\Omega)$ to a limit denoted by $Z_\infty$, which belongs to the $q$th chaos of $B^H$.
\end{prop}
{\it Proof}.
Let us first prove the convergence in $L^2(\Omega)$.
For $n,m\geq 1$, we have
\[
E[Z_nZ_m]=q!(nm)^{q-1}\sum_{k=0}^{n-1}\sum_{l=0}^{m-1}\left(E\big[\big(B^H_{(k+1)/n}-B^H_{k/n}\big)
\big(B^H_{(l+1)/m}-B^H_{l/m}\big)\big]\right)^q.
\]
Furthermore, since $H>1/2$, we have for all $s,t\geq 0$,
\[
E[B^H_sB^H_t]=H(2H-1)\int_0^t du\int_0^s dv |u-v|^{2H-2}.
\]
Hence
\[
E[Z_nZ_m]=q!H^q(2H-1)^q\times\frac1{nm}
\sum_{k=0}^{n-1}\sum_{l=0}^{m-1}\left(nm\int_{k/n}^{(k+1)/n} du \int_{l/m}^{(l+1)/m} dv
|v-u|^{2H-2}\right)^q.
\]
Therefore, as $n,m\to\infty$, we have,
\[
E[Z_nZ_m]\to
q!H^q(2H-1)^q
\int_{[0,1]^2}|u-v|^{(2H-2)q}dudv,
\]
and the limit is finite since $H>1-\frac{1}{2q}$.
In other words, the sequence $\{Z_n\}$ is Cauchy in $L^2(\Omega)$,
and hence converges in $L^2(\Omega)$ to some $Z_\infty$.

Let us now prove that $\{Z_n\}$ converges also almost surely.
Observe first that, since $Z_n$ belongs to the $q$th chaos of $B^H$ for all $n$, since $\{Z_n\}$
converges in $L^2(\Omega)$ to $Z_\infty$ and since the $q$th chaos of $B^H$ is closed in $L^2(\Omega)$ by definition, we have that $Z_\infty$ also belongs to the
$q$th chaos of $B^H$.
In \cite[Proposition 3.1]{BN}, it is shown that $E[|Z_n-Z_\infty|^2]\leq Cn^{2q-1-2qH}$, for some
positive constant $C$ not depending on $n$.
Inside a fixed chaos, all the $L^p$-norms are equivalent. Hence, for any $p>2$,
we have $E[|Z_n-Z_\infty|^p]\leq Cn^{p(q-1/2-qH)}$.
Since $H>1-\frac{1}{2q}$, there exists $p>2$ large enough such that
$(q-1/2-qH)p<-1$. Consequently
\[
\sum_{n\geq 1}E[|Z_n-Z_\infty|^p]<\infty,
\]
leading, for all $\e>0$, to
\[
\sum_{n\geq 1}P[|Z_n-Z_\infty|>\e]<\infty.
\]
Therefore, we deduce from the Borel-Cantelli lemma
that $\{Z_n\}$ converges
almost surely to $Z_\infty$.
\fin

We now want to see if one can associate almost sure central limit theorems to the convergences in law
(\ref{eq:Breuer_Major1}), (\ref{eq:Breuer_Major2}) and (\ref{eq:Breuer_Major3}).
We first consider the case $H<1-\frac{1}{2q}$.
\begin{prop}\label{BMsouscritique}
Assume that $q\geq 2$ and that $H< 1-\frac{1}{2q}$, and consider
\[
G_n=\frac{V_n}{\sigma_n\,\sqrt{n}}
\] as in (\ref{eq:Breuer_Major1}).
Then, $\{G_n\}$ satisfies an ASCLT.
\end{prop}
{\it Proof}.
Since $2H-2>1/q$, it suffices to combine (\ref{behaviorrho}), Proposition \ref{prop-az}
and Theorem \ref{BM}.
\fin

Next, let us consider the critical case $H=1-\frac1{2q}$.
In this case,
$
\sum_{r\in\mathbb{Z}} |\rho(r)|^q=\infty.
$
Consequently, as it is impossible to apply Theorem \ref{BM},
we propose another strategy which relies on
the following lemma established in \cite{BN}.
\begin{lemma}\label{Stein-lemma}
Set $H=1-\frac{1}{2q}$. Let $\HH$ be the real and separable Hilbert space
defined as follows: (i) denote by $\mathscr{E}$ the
set of all $\mathbb{R}$-valued step functions on $[0,\infty)$, (ii)
define $\EuFrak H$ as the Hilbert space obtained by closing
$\mathscr{E}$ with respect to the scalar product
\[
\left\langle
{\mathbf{1}}_{[0,t]},{\mathbf{1}}_{[0,s]}\right\rangle _{\EuFrak
H}=E[B^H_tB^H_s].
\]
For any $n\geq 2$, let $f_n$ be the element of $\HH^{\odot q}$ defined by
\begin{equation}\label{fk}
f_n=
\frac{1}{\sigma_n\,\sqrt{n\log n}}\sum_{k=0}^{n-1}{\bf 1}_{[k,k+1]}^{\otimes q},
\end{equation}
where $\sigma_n$ is the positive normalizing constant which ensures that
$q!\|f_n\|^2_{\HH^{\otimes q}}=1$.
Then, there exists a constant $C>0$, depending only on $q$ and $H$  such that,
for all $n\geq 1$ and $r=1,\ldots,q-1$
\[
\|f_n\otimes_r f_n\|_{\HH^{\otimes (2q-2r)}}\leq C
(\log n)^{-1/2}.
\]
\end{lemma}
We can now state and prove the following result.
\begin{prop}\label{BMcritique}
Assume that $q\geq 2$ and $H= 1-\frac{1}{2q}$, and consider \[
G_n=\frac{V_n}{\sigma_n\,\sqrt{n\log n}}
\] as in (\ref{eq:Breuer_Major2}).
Then, $\{G_n\}$ satisfies an ASCLT.
\end{prop}
{\bf Proof of Proposition \ref{BMcritique}}. We shall make use of Corollary \ref{main}.
Let $C$ be a positive constant, depending only on $q$ and $H$, whose value may change from line to line.
We consider the real and separable Hilbert space $\EuFrak H$ as defined in Lemma \ref{Stein-lemma}.
We have $G_n=I_q(f_n)$ with $f_n$ given by (\ref{fk}). According to Lemma \ref{Stein-lemma}, we have
for all $k\geq 1$ and $r=1,\ldots,q-1$, that
$\|f_k\otimes_r f_k\|_{\HH^{\otimes (2q-2r)}}\leq C(\log k)^{-1/2}$.
Hence
\begin{eqnarray*}
\sum_{n\geq 2}\frac{1}{n\log^2 n}\sum_{k=1}^n \frac{1}{k}\|f_k\otimes_r f_k\|_{\HH^{\otimes (2q-2r)}}
&\leq&
C\sum_{n\geq 2}\frac{1}{n\log^2 n}\sum_{k=1}^n \frac{1}{k\sqrt{\log k}},\\
&\leq&
C\sum_{n\geq 2}\frac{1}{n\log^{3/2} n}
<\infty.
\end{eqnarray*}
Consequently, assumption $(A'_1)$ is satisfied.
Concerning $(A'_2)$, note that
\[
\langle f_k,f_l\rangle_{\HH^{\otimes q}} =\frac{1}{\sigma_k\sigma_l\,\sqrt{k\log k}\sqrt{ l\log l}}
\sum_{i=0}^{k-1}\sum_{j=0}^{l-1}\rho(j-i)^q.
\]
We deduce from Lemma \ref{sigma} below that $\sigma_n^2\rightarrow\sigma_\infty^2>0$.
Hence, for all $l\geq k\geq 1$
\begin{eqnarray*}
\big|\langle f_k,f_l\rangle_{\HH^{\otimes q}}\big|&\leq&\frac{C}{\sqrt{k\log k}\sqrt{l\log l}}
\sum_{i=0}^{k-1}\sum_{j=0}^{l-1}\big|\rho(j-i)\big|^q,\\
&=&
\frac{C}{\sqrt{k\log k}\sqrt{l\log l}}
\sum_{i=0}^{k-1}\sum_{r=-i}^{l-1-i}\big|\rho(r)\big|^q,\\
&\leq&C\frac{\sqrt{k}}{\sqrt{\log k}\sqrt{l\log l}}\,\sum_{r=-l}^{l}\big|\rho(r)\big|^q
\leq C\sqrt{\frac{k\log l}{l\log k}}.
\end{eqnarray*}
The last inequality follows from the fact that $\sum_{r=-l}^l\big|\rho(r)\big|^q\leq C\log l$
since, by (\ref{behaviorrho}), as $|r|\to\infty$,
\[
\rho(r)\sim\left(1-\frac{1}{q}\right)\left(1-\frac1{2q}\right)|r|^{-1/q}.
\]
Finally, assumption $(A'_2)$ is also satisfied as
\begin{eqnarray*}
\sum_{n\geq 2}\frac{1}{n\log^3 n}\sum_{k,l=2}^{n} \frac{\big|\langle f_k,f_l\rangle_{\HH^{\otimes q}}\big|}{kl}
&\leq&
2\sum_{n\geq 2}\frac{1}{n\log^3 n}\sum_{l=2}^{n}\sum_{k=2}^l \frac{\big|\langle f_k,f_l\rangle_{\HH^{\otimes q}}\big|}{kl},\\
&\leq&C\sum_{n\geq 2}\frac{1}{n\log^3 n}\sum_{l=2}^n \frac{\sqrt{\log l}}{l^{3/2}}\sum_{k=2}^l \frac1{\sqrt{k\log k}},\\
&\leq&C\sum_{n\geq 2}\frac{1}{n\log^3 n}\sum_{l=2}^n \frac1l   \leq
C\sum_{n\geq 2}\frac{1}{n\log^2 n}<\infty.
\end{eqnarray*}
\fin

In the previous proof, we used the following lemma.
\begin{lemma}\label{sigma}
Assume that $q\geq 2$ and $H=1-\frac1{2q}$.
Then,
\[
\sigma_n^2\to 2q!\left(1-\frac1q\right)^q\left(1-\frac1{2q}\right)^q >0,\quad\mbox{as $n\to\infty$.}
\]
\end{lemma}
{\bf Proof}.
We have $E[(B^H_{k+1}-B^H_k)(B^H_{l+1}-B^H_l)]=\rho(k-l)$
where $\rho$ is given in (\ref{e:cov}).
Hence,
\begin{eqnarray*}
E[V_n^2]&=&\sum_{k,l=0}^{n-1}E\big(H_q(B^H_{k+1}-B^H_k)H_q(B^H_{l+1}-B^H_l)\big)
=q!\sum_{k,l=0}^{n-1}\rho(k-l)^q,\\
&=&q!\sum_{l=0}^{n-1}\sum_{r=-l}^{n-1-l}\rho(r)^q
=q!\sum_{|r|<n} \big(n-1-|r|\big)\rho(r)^q,\\
&=&q!\left(n\sum_{|r|<n} \rho(r)^q - \sum_{|r|<n} \big(|r|+1\big)
\rho(r)^q\right).
\end{eqnarray*}
On the other hand, as $|r|\to\infty$,
\[
\rho(r)^q\sim
\left(1-\frac1q\right)^q\left(1-\frac1{2q}\right)^q\frac{1}{|r|}.
\]
Therefore, as $n\to\infty$,
\[
\sum_{|r|<n} \rho(r)^q \sim
\left(1-\frac1{2q}\right)^q\left(1-\frac1q\right)^q\sum_{0<|r|<n}\frac{1}{|r|}
\sim 2\left(1-\frac1{2q}\right)^q\left(1-\frac1q\right)^q \log n
\]
and
\[
\sum_{|r|<n} \big(|r|+1\big)\rho(r)^q \sim
\left(1-\frac1{2q}\right)^q\left(1-\frac1q\right)^q\sum_{|r|<n} 1
\sim 2n
\left(1-\frac1{2q}\right)^q\left(1-\frac1q\right)^q.
\]
Consequently, as $n\to\infty$,
\[
\sigma_n^2=\frac{E[V_n^2]}{n \log n}\to 2q!\left(1-\frac1q\right)^q\left(1-\frac1{2q}\right)^q.
\]
\fin

Finally, we consider
\begin{equation}\label{e:c3}
G_n=n^{q(1-H)-1}V_n
\end{equation}
with $H>1-\frac{1}{2q}$. We face in this case some difficulties. First,
since the limit of $\{G_n\}$ in (\ref{eq:Breuer_Major3}) is not Gaussian, we cannot apply our general
criterion Corollary \ref{main} to obtain an ASCLT.
To modify adequably the criterion, we would need a version of Lemma \ref{noupec}
for random variables with an Hermite distribution, a result which is not presently available.
Thus, an ASCLT associated to the convergence in law (\ref{eq:Breuer_Major3}) falls outside the
scope of this paper. We can nevertheless make a number of observations.
First, changing the nature of the random variables without changing their law has no impact
on CLTs as in (\ref{eq:Breuer_Major3}), but may have a great impact on an ASCLT. To see this,
observe that for each fixed $n$, the ASCLT involves not only the distribution of the single variable $G_n$, but also
the joint distribution of the vector $(G_1,\ldots,G_n)$.

Consider, moreover, the following example. Let $\{G_n\}$ be a sequence of random variables converging in law
to a limit $G_\infty$. According to a theorem of Skorohod, there is a sequence $\{G_n^*\}$ such that
for any fixed $n$, $G_n^*\overset{\rm law}{=}G_n$ and such that $\{G_n^*\}$ converges almost surely,
as $n\to\infty$, to a random variable $G^*_\infty$ with $G_\infty^*\overset{\rm law}{=}G_\infty$.
Then, for any bounded continuous function $\varphi:\R\to\R$, we have
$\varphi(G_n^*)\longrightarrow\varphi(G_\infty^*)$ a.s. which clearly implies the
almost sure convergence
\[
\frac1{\log n}\sum_{k=1}^n \frac1k \varphi(G_k^*)\longrightarrow\varphi(G_\infty^*).
\]
This limit is, in general, different from $E[\varphi(G_\infty^*)]$ or equivalently
$E[\varphi(G_\infty)]$,
that is, different from the limit if one had an ASCLT.

Consider now the sequence $\{G_n\}$ defined by (\ref{e:c3}).

\begin{prop}
The Skorohod
version of
\begin{equation}\label{e:Gn}
G_n=n^{q(1-H)-1}\sum_{k=0}^{n-1}H_q(B^H_{k+1}-B^H_k)
\end{equation}
is
\begin{equation}\label{e:Gnstar}
G_n^*=Z_n=n^{q(1-H)-1}\sum_{k=0}^{n-1}H_q\big(
n^H(B^H_{(k+1)/n}-B^H_{k/n})\big),
\end{equation}
\end{prop}
{\it Proof}. Just observe that $G_n^*\overset{\rm law}{=}G_n$ and
$G_n^*$ converges almost surely by Proposition \ref{lm-hermite-yahoo}.
\fin

Hence, in the case of Hermite distributions, by suitably modifying the argument of the Hermite polynomial
$H_q$ in a way which does not change the limit in law, namely by considering $Z_n$ in (\ref{e:Gnstar})
instead of $G_n$ in (\ref{e:Gn}), we obtain the almost sure convergence
\[
\frac1{\log n}\sum_{k=1}^n\frac1k \varphi(Z_k)\longrightarrow
\varphi(Z_\infty).
\]
The limit $\varphi(Z_\infty)$ is, in general, different from the limit
expected under an ASCLT, namely $E[\varphi(Z_\infty)]$, because $Z_\infty$ is a non-constant
random variable with an Hermite distribution
(Dobrushin and Major \cite{DoMa}, Taqqu \cite{T}).
 Thus, knowing the law of $G_n$ in (\ref{e:Gn}), for a fixed $n$,
does not allow to determine whether an ASCLT holds or not.
\\

\bigskip

\noindent
{\bf Acknowledgments}. This paper originates from the conference ``Limit theorems and applications'',
University Paris I Panth\'eon-Sorbonne,
January 14-16, 2008, that the three authors were attending. We warmly thank J.-M. Bardet and C. A. Tudor
for their invitation and generous support. Also,
I. Nourdin would like to thank M. S. Taqqu
for his hospitality during his stay at Boston University in March 2009, where part of this research was carried out.
Finally, we would like to thank two anonymous referees for their careful reading of the manuscript and for their
valuable suggestions and remarks.

\end{document}